%Paper: 9205004
%From: Takao Fujita <fujita@math.titech.ac.jp>
%Date: Fri, 8 May 92 14:59:24 JST

\input amstex
\documentstyle{amsppt}

\topmatter
\title Notes on Kodaira energies of polarized varieties
\endtitle
\author Takao FUJITA \endauthor
\address {Department of Mathematics
\newline
Tokyo Institute of Technology
\newline
Oh-okayama, Meguro, Tokyo
\newline
152 Japan}
\endaddress
\endtopmatter

\define\dnl{\newline\newline}
\define\lra{{\longrightarrow}}
\define\nl{\newline}

\define\Inf{{\roman{Inf}}}

%\define\det{{\roman{det}}}

\define\dm{{\roman{dim}}}

\define\SE{\Cal E}

\define\SO{\Cal O}

\define\BC{{\Bbb C}}
\define\BP{{\Bbb P}}
\define\BQ{{\Bbb Q}}
\define\BZ{{\Bbb Z}}

\define\Rmk{{\it Remark. }}

\define\Lim{{\roman{Lim}}}
\define\ke{\kappa\epsilon}
\define\fl{^\flat}

\document
In this note we propose a couple of conjectures concerning Kodaira energies of
polarized varieties and give a few partial answers.\dnl
{\bf \S1. Conjectures}

Let $V$ be a variety over $\BC$ and let $B=\sum b_iB_i$ be an effective
$\BQ$-Weil divisor on $V$ such that $b_i\le1$ for any $i$.
Such a pair $(V, B)$ will be called a {\it log variety}.
It is said to be {\it log terminal} if it has only weak log terminal
singularities in the sense of [{\bf KMM}].
In this case, the $\BQ$-bundle $K_V+B$ is called the log canonical bundle of
$(V, B)$ and will be denoted by $K(V,B)$.

A $\BQ$-bundle $L$ on a log terminal variety $(V, B)$ is said to be {\it log
ample} if there is an effective $\BQ$-divisor $E$ such that $(V, B+E)$ is log
terminal and $L-\epsilon E$ is ample for any $0<\epsilon\le 1$.
Note that ``log ample'' implies ``nef big'', and the converse is also true if
$b_i<1$ for all $i$.

For a big $\BQ$-bundle $L$ on a log terminal variety $(V, B)$, we define
$$\kappa\epsilon(V, B, L)=-\Inf\{t\in\BQ\vert\kappa(K(V,B)+tL)\ge 0\},$$
which will be called the {\it Kodaira energy} of $(V, B, L)$.
When $B=0$, we write simply $\kappa\epsilon(V, L)$.

Clearly $\kappa\epsilon(V, B, L)<0$ if and only if $K(V,B)$ is not
pseudo-effective.
We conjecture that $\kappa\epsilon(V, B, L)\in\BQ$ in this case (cf. [{\bf
Ba}]).
This will be derived from the following

{\bf Fibration Conjecture.} {\sl Let $L$ be a big $\BQ$-bundle on a log
terminal variety $(V,B)$ such that $k=\ke(V,B,L)<0$ and $K(V,B)+tL$ is log
ample for many $t>0$.
Then there is a birational model $(V', B')$ of $(V, B)$ together with a
morphism $\Phi: V'\lra W$ such that $\dm W<\dm V$, $\Phi_*\SO_{V'}=\SO_W$, the
relative Picard number $\rho(V'/W)$ is $1$, and $K(V',B')-kL'=\Phi^*A$ for some
ample $\BQ$-bundle $A$ on $W$, where $L'$ is the proper transform (as a Weil
divisor) of $L$ on $V'$.}

The birational map $V\lra V'$ will be obtained by applying the Log Minimal
Model Programm, and will be a composite of elementary divisorial contractions
and flips.
Note that $L'$ may not be nef even if $L$ is ample.

Any way, we have $K(F, B'_F)=kL'_F$ for any general fiber $F$ of $\Phi$, and
$L'_F$ is ample since $\rho(V'/W)=1$.
Hence, by a certain conjectural boundedness of $\BQ$-Fano varieties, we shall
obtain the following

{\bf Spectrum Conjecture.} {\sl Let $S_n$ be the Kodaira spectrum of polarized
$n$-folds, namely, the set of all the possible Kodaira energies of $(V, L)$,
where $V$ is a variety with $\dm V=n$ having only terminal singularities and
$L$ is an ample line bundle on $V$.
Then $\{t\in S_n\vert t<-\delta\}$ is a finite subset of $\BQ$ for any
$\delta>0$.}

If we allow $V$ to have log terminal singularities, the assertion of the
Spectrum Conjecture is false.

To be precise, let $\Lim(X)$ denote the set of limit points of $X$, namely,
$p\in\Lim(X)$ if and only if $U\cap X$ is an infinite set for any neighborhood
$U$ of $p$.
Let $\Lim{}^k(X)=\Lim(\Lim{}^{k-1}(X))$, let $X\sqcup Y$ denote
$(X-Y)\cup(Y-X)=(X\cup Y)-(X\cap Y)$ and let $S'_n$ be the set of all the
possible Kodaira energies of $(V, B, L)$ such that $(V, B)$ is log terminal,
$\dm V=n$, $B$ is a $\BZ$-Weil divisor (or equivalently, $b_i=1$ for any $i$;
possibly $B=0$), and $L$ is ample on $V$.
Then we have the following

{\bf Log Spectrum Conjecture.} {\sl For any $k\le n$, let $S'_{n,k}$ be the set
$\{t\vert t-k\in S'_n {\text{ and }} -1<t\le0\}$.
Then $\Lim{}^{n-k}(S'_{n,k})=\{0\}$ and $\Lim{}^{n-k}(S'_{n,k}\sqcup
S'_{n+1,k+1})=\Lim{}^{n-k}(S'_{n,k}\sqcup \Lim(S'_{n,k-1}))=\emptyset$.
Moreover, for any $s<0$, there exists $\delta>0$ such that $\{t\in S'_n\vert
s<t<s+\delta\}$ is a finite set.}

The conclusion cannot be simplified even if $(V, B)$ is assumed to be smooth,
or if we assume $B=0$ allowing $V$ to have log terminal singularities.
\dnl
{\bf \S2. Results}

The preceding conjectures are verified if $\dm V\le 3$ under mild additional
assumptions.

{\bf Theorem 1.} {\sl Let $(V, B, L)$ be as in the Fibration Conjecture.
Suppose that $n=\dm V\le 3$ and that $V$ is $\BQ$-factorial, namely, every Weil
divisor on $V$ is $\BQ$-Cartier. Then the assertion of the conjecture is true.}

{\bf Theorem 2.} {\sl The Spectrum Conjecture is true for the Kodaira spectrum
of polarized 3-folds such that $V$ is $\BQ$-factorial.}

\Rmk The $\BQ$-factoriality is needed to apply the theory [{\bf Sho}], [{\bf
Ka2}], [{\bf Ka1}].

{\it Outline of proof of Theorem 1.}
Set $\tau=\Inf\{t\in\BQ\vert K(V,B)+tL \text{ is nef}\}$.
Then we see $\tau\in\BQ$ by using Cone Theorem (cf. [{\bf KMM}]).
Next, by using the Base-Point-Free Theorem, we get a fibration $f: V\lra X$ and
an ample $\BQ$-bundle $A$ on $X$ such that $K(V,B)+\tau L=f^*A$.

Now we let the Log Minimal Model Programm (cf. [{\bf KMM}], [{\bf Sho}], [{\bf
Ka2}]) run over $X$.
By several elementary divisorial contractions and log flips, $(V, B)$ is
transformed to a pair $(V_1, B_1)$ satisfying one of the following
conditions:\nl
(1) There is an extremal ray $R$ on $V_1$ over $X$ such that its contraction
morphism $\rho: V_1\lra W$ is of fibration type.\nl
(2) $K(V_1, B_1)$ is relatively nef over $X$.

During the process, $L$ is transformed to a $\BQ$-bundle $L_1$ on $V_1$ as Weil
divisors.
It is easy to see that the bigness is preserved, and that the Kodaira energy
does not change.
Thus, in case (1), we are done by setting $V'=V_1$.

In case (2), we can show that $A_{V_1}$ is log ample on $(V_1, B_1)$.
This is easy to prove when $b_i<1$ for all $i$, but the proof is a little
complicated in general.

Thus, $(V_1, B_1, L_1)$ satisfies the same condition as $(V, B, L)$.
Note also that $\tau_2=\Inf\{t\in\BQ\vert K(V_1,B_1)+tL_1 \text{ is
nef}\}<\tau$.
By the same process as above we get another triple $(V_2, B_2, L_2)$, and
continue as long as necessary.
By the termination theorem, we reach the above situation (1) after finite
steps.

For the proof of Theorem 2, [{\bf Ka1}] is essential.
\dnl
{\bf \S3. Classification}

By the same method as in [{\bf F2}], we can classify smooth polarized 3-folds
$(M, L)$ with $\ke(M, L)<-\frac12$ as follows.

(3.1) $K+3L$ is nef and $\ke(M, L)\ge -3$ unless $(M, L)\cong(\BP^3, \SO(1))$.

(3.2) $K+2L$ is nef and $\ke\ge -2$ unless $(M, L)$ is a smooth scroll over a
curve or a hyperquadric in $\BP^4$.
$\ke=-3$ in these cases.

(3.3) From now on, $K+2L$ is assumed to be nef.
If there is a divisor $E$ such that $(E, L_E)\cong(\BP^2, \SO(1))$ and
$[E]_E=\SO(-1)$, let $\pi: M\lra M_1$ be the blow down of $E$ to a smooth
point.
Then the push-down $L_1$ of $L$ is ample on $M_1$ and $\pi^*(K_1+2L_1)=K+2L$
for the canonical bundle $K_1$ of $M_1$.
If there is a similar divisor on $M_1$, we blow it down again.
After several steps, we get a polarized manifold $(M', L')$ such that
$(K'+2L')_M=K+2L$, $\ke(M', L')=\ke(M, L)\ge -2$ on which there is no divisor
of the above type.
This model $(M', L')$ is called the (first) reduction of $(M, L)$.

(3.4) $\ke(M, L)=-2$ if and only if $K+2L$ is not big.
In this case, according to the value of $\kappa(K+2L)$, $(M, L)$ is classified
as follows:\nl
(3.4.0) $K+2L=0$, i.e., $(M, L)$ is a Del Pezzo 3-fold.\nl
(3.4.1) $(M, L)$ is a hyperquadric fibration over a curve.\nl
(3.4.2) $(M, L)$ is a scroll over a surface.

\Rmk In the cases (3.4.1) and (3.4.2), we have $M=M'$ by the ampleness of $L$.

(3.5) $K+2L$ is nef and big if and only if $K'+2L'$ is ample.
Moreover, in this case, $K'+L'$ is nef except the following cases:\nl
(3.5.1) $M'$ is a $\BP^2$-bundle over a curve and $L'_F=\SO(2)$ for any fiber
$F$. $\ke=-3/2$.\nl
(3.5.2) $M'$ is a hyperquadric in $\BP^4$ and $L'=\SO(2)$. $\ke=-3/2$.\nl
(3.5.3) $(M', L')\cong(\BP^3, \SO(3))$. $\ke=-4/3$.

(3.6) From now on, $K'+L'$ is assumed to be nef.
Then it is not big if and only if $\ke(M, L)=-1$.
These cases are classified as follows:\nl
(3.6.0) $K'+L'=0$.\nl
(3.6.1) $(M', L')$ is a Del Pezzo fibration over a curve.\nl
(3.6.2) $(M', L')$ is a conic bundle over a surface.

(3.7) To study the case in which $K'+L'$ is nef big, we use the theory of
second reduction as in [{\bf BS}].
We have a birational morphism $\Phi: M'\lra M''$ such that $K'+L'=\Phi^*A$ for
some ample line bundle $A$ on $M''$, and this pair $(M'', A)$ is called the
second reduction of $(M, L)$.
However, unlike the case of first reduction, $M''$ may have singularities,
$L''=\Phi_*L'$ may not be invertible, may not be nef.
By a careful analysis of the map $\Phi$ using Mori theory, we see that the
singularity of $M''$ is of very special type.
It is a hypersurface singularity of the type $\{x^2+y^2+z^2+u^k=0\}\ (k=2, 3)$,
or the quotient singularity isomorphic to the vertex of the cone over the
Veronese surface $(\BP^2, \SO(2))$.
In particular $L''$ is invertible except at quotient singularities and $2L''$
is invertible everywhere.

The cases $-1<\ke<-1/2$ can be classified according to the type of the second
reduction $(M'', A)$ as follows.

(3.8) $\ke=-4/5$. $(M'', A)\cong(\BP^3, \SO(1))$ and $L''=\SO(5)$.

(3.9.0) $\ke=-3/4$. $M''$ is a hyperquadric in $\BP^4$ and $L''=\SO(4)$.

(3.9.1) $\ke=-3/4$. $(M'', A)$ is a scroll over a curve. $L''_F=\SO(4)$ for any
fiber $F$.

(3.10) $\ke=-5/7$. $(M'', A)$ is a cone over $(\BP^2, \SO(2))$.
Compare [{\bf F2};(4.8.0)].

(3.11.0) $\ke=-2/3$. $(M'', A)$ is a Del Pezzo 3-fold, i.e., $K''=-2A$.

(3.11.1) $\ke=-2/3$. $(M'', A)$ is a hyperquadric fibration over a curve.

(3.11.2) $\ke=-2/3$. $(M'', A)$ is a scroll over a surface.

In the following cases, $(M'', A)$ may need to be further blow down to another
model $(M^\flat, A\fl)$.
But this pair has no worse singularities than $(M'', A)$.

(3.12.0.0) $\ke=-3/5$. $M\fl$ is a hyperquadric in $\BP^4$ and $A\fl=\SO(2)$.

(3.12.0.1) $\ke=-3/5$. $M\fl$ has exactly one quotient singularity, and the
blow-up $M^\#$ at this point is isomorphic to the blow-up of $\BP^3$ along a
smooth plane cubic $C$.
$A\fl_{M^\#}$ is $3H-E_C$, where $H$ is the pull-back of $\SO(1)$ of $\BP^3$
and $E_C$ is the exceptional divisor over $C$.
Compare [{\bf F2};(4.6.0.1.0)].

(3.12.0.2) $\ke=-3/5$. $M\fl$ has exactly two quotient singularities, and the
blow-up $M^\#$ at these points is a smooth member of $\vert 2H(\SE)\vert$ on
$\BP(\SE)$, where $\SE$ is the vector bundle $\SO(2)\oplus\SO\oplus\SO$ on
$\BP^2$ and $H(\SE)$ is the tautological line bundle on $\BP(\SE)$.
$A\fl_{M^\#}$ is the restriction of $H(\SE)$.
Compare [{\bf F2};(4.6.0.2.1)].

(3.12.1) $\ke=-3/5$. $M\fl$ is a $\BP^2$-fibration over a curve and
$A\fl_F=\SO(2)$ for any general fiber $F$.

(3.13) $\ke=-4/7$. $(M\fl, A\fl)\cong(\BP^3, \SO(3))$ and $L\fl=\SO(7)$.

(3.14) $\ke=-5/9$. $(M\fl, B)$ is a cone over $(\BP^2, \SO(2))$ for some line
bundle $B$ such that $A\fl=2B$.

(3.15) When $\ke\ge -1/2$, we can reduce the problem to the case in which
$K\fl+A\fl$ is nef, where $K\fl$ is the canonical $\BQ$-bundle of $M\fl$.

Note that $M\fl$ is obtained from $M$ without using flips.

\enddocument